# Diffusion bridge with misspecification: theory construction and application to high-resolution fish count data

Hidekazu Yoshioka[1, *]

[1]Japan Advanced Institute of Science and Technology, 1-1 Asahidai, Nomi, Ishikawa 923-1292, Japan. ORCID: 0000-0002-5293-3246

*Corresponding author: yoshih@jaist.ac.jp, +81-761-51-1745

**Abstract**

Stochastic processes of bridge types having pinned initial and terminal conditions have been widely used in applied research areas, but they all have a common drawback in that the model at hand is possibly misspecified owing to its stochastic nature; namely, parameter values and coefficients are distorted compared to the ground truth. We consider a pair of novel exactly-solvable optimization problems that provide both the lower and upper bounds of the performance index of a diffusion bridge. Our formulation is based on the Girsanov transformation, in which the model uncertainty is measured through relative entropy. We provide a sufficient condition under which these optimization problems are well-posed, and hence admit the corresponding maximizer/minimizer that achieves the worst-case lower and upper bounds given the ambiguity aversion or uncertainty size. We apply the proposed method to the latest 10-min, high-resolution fish count data of a migratory fish in a river and discuss the influence of model uncertainty on the estimation of the total fish count, which is an important problem in resource and environmental management.

***Statements & declarations***

**Acknowledgments:** The author thanks the Japan Water Agency, Ibi River, and Nagara River General Management Offices for providing valuable fish migration data from recent years.

**Funding:** This work was supported by the Japan Science and Technology Agency (PRESTO No. JPMJPR24KE).

**Conflict of interests:** The author declares no conflicts of interest.

**Data availability:** The data will be made available upon reasonable request to the corresponding author.

**Declaration of generative AI in scientific writing:** The author did not use generative AI in the scientific writing of the manuscript.

**Contribution:** The author prepared all parts of this manuscript.

## 1. Introduction

### 1.1 Background

Stochastic differential equations (SDEs) are ordinary differential equations (ODEs) driven by martingale noises [1], such as Brownian motions and Poisson processes, and have been applied to a wide variety of phenomena in natural sciences, ranging from cell molecular chemistry [2], system biology [3], population dynamics [4], and fluid dynamics [5] to cosmology [6]. SDEs are also frequently used in social sciences, such as finance [7], insurance [8], and opinion dynamics [9]. Machine learning is another field of research in which SDEs play a vital role [10,11].

SDEs with pinned initial and terminal conditions are referred to as bridges, particularly those driven by Brownian motion (i.e., time-continuous noise) as diffusion bridges. The simplest diffusion bridge is the Brownian bridge, which is a Brownian motion with pinned terminal conditions [12]. More sophisticated models also exist where underlying stochastic processes are governed by linear [13,14], fractional [15], or even more strongly nonlinear SDEs [16,17] depending on the target problem. The application of bridges to efficient sampling problems has also been investigated [18,19]. Diffusion bridges are fundamentally different from classical SDEs because the latter do not necessarily pin the terminal conditions. A fundamental theoretical difficulty in a diffusion bridge arises because of the pinning of a terminal condition, as the behavior of its solutions is constrained at and near the terminal time. Typically, pinning a terminal condition in a diffusion bridge is achieved by suitably specifying exploding drift and diffusion coefficients in the corresponding SDE so that the terminal condition is continuously reached with probability one [20]. Therefore, a balance between the drift and diffusion coefficients is required when modeling a diffusion bridge.

A topic that has not been well-studied concerning bridges is model misspecification, that is, studying the performance of a model under the assumption that its parameter values and/or functional shapes of coefficients may differ from the truth [21,22,23,24]. Misspecified SDEs provide quantitatively and often qualitatively different outputs compared to the ground truth. This issue is serious for SDEs in both theory and application because they are often efficient but conceptual, simplified mathematical models of complex phenomena. Operating misspecified models along with their uncertainties is therefore a key topic in modeling with SDEs. From a practical standpoint, such an approach should not lose the efficiency of SDEs; moreover, it should be mathematically rigorous.

The fundamental theory for modeling with possibly misspecified SDEs was established by Hansen and Sargent [25], who considered the stochastic control problems of misspecified SDEs, where the degree of their uncertainties was measured using the relative entropy between true and misspecified models. The key here is that the benchmark and misspecified models are connected by a specific transformation of probability measures, called Girsanov transformation in the field of stochastic calculus (e.g., Chapter 1.4 in Øksendal and Sulem [1]), which arises as an additional drift coefficient in a misspecified SDE. An advantage of this framework is that solving a control problem with a misspecified model reduces to finding a proper solution to an optimality equation, as in classical stochastic control problems; therefore, existing analytical and computational techniques for solving stochastic control problems (e.g., Chapter 5 in Øksendal and Sulem [1]; Chapter 3 in Pham [26]) often carry over to cases with misspecifications. With this remarkable property, the framework of Hansen and Sargent [25] and related studies have been successfully employed in a broad range of research fields, such as economics [27,28,29], climate change [30,31], and environmental conservation [32]. However, their approach for modeling with misspecified SDEs has not been well studied, except for Baltas and Yannacopoulos [33], who studied the Brownian bridge with misspecification, focusing on its application to informed traders, although diffusion bridges have been employed in a variety of situations. This motivated us to conduct this study with the aims and contributions explained below.

### 1.2 Aim and contribution

The aim of this study is to formulate and theoretically analyze a diffusion bridge with possible misspecifications. We also address the engineering problem so that this study seamlessly covers its theory and application. In short, the main objective of this study is to formulate and analyze optimization problems that give upper and lower bounds of the expectation of total fish count, given the uncertainty in the coefficients of SDE $(1)$ given below. More specifically, we show that solving these optimization problems reduces to the study of certain ODEs. We also investigate the sufficient conditions for the optimality of these ODEs. The contribution of this study to achieving its aims is explained in this subsection.

The diffusion bridge we study is the Cox–Ingersoll–Ross (CIR) bridge presented below (formulated rigorously with proper notations in the next section).

$$\mathrm{d}X_t = \left( a_t - \frac{r_t}{1-t} X_t \right) \mathrm{d}t + \sigma_t \sqrt{\frac{r_t}{1-t} X_t} \mathrm{d}B_t , \quad 0 < t < 1 , \quad X_0 = X_1 = 1 , \tag{1}$$

where $t$ is time, $B = (B_t)_{t \geq 0}$ is a 1-D standard Brownian motion, the diffusion term is defined in Itô's sense, and $a, r, \sigma > 0$ are some functions of time. This scalar diffusion bridge with homogeneous initial and terminal conditions was proposed by Yoshioka [34] to model the unit-time count data of migrating fish along a river and admits a pathwise unique solution $X = (X_t)_{0 \leq t \leq 1}$ that is nonnegative; however, its misspecified versions have not been studied thus far. This bridge, called the CIR bridge, has drift and diffusion coefficients that explode at the terminal time $t = 1$, but its solution attains the terminal condition continuously with a probability of one. Balancing the explosion speed in the drift and diffusion coefficients in $(1)$ is the key to this property.

A fundamental quantity to be modeled in an industrial application of the CIR bridge $(1)$ is the total fish count during a certain time interval, for example, one day, which is a key statistic for assessing fish populations [35,36,37,38,39]. Identification of an SDE in fisheries and ecological applications often faces modeling errors, which are considered a form of misspecification owing to the lack of quantity and/or quality of data [e.g., 40,41,42]. In practice, population assessment of migrating fish, both over- and underestimations, would be possible; hence, their estimations, given a degree of uncertainty, would be helpful in applications such as fisheries science.

The principle in the mathematical framework in this study is that the uncertainty of the model stems from distortions in the Brownian motion, and this assumption contributes to the high analytical tractability of the proposed framework. In this study, using the framework of Hansen and Sargent [25], we formulate a pair of optimization problems to obtain upper- and lower-bound estimations of the population, given the size of ambiguity aversion or relative entropy, which models the difference between the benchmark and misspecified models (**Figure 1** is an overview diagram of the proposed framework). Each optimization problem, which is a minimization or maximization problem of an expectation subject to an SDE, is reduced to finding an appropriate solution to an optimality equation, which is a nonlinear degenerate parabolic partial differential equation with exploding coefficients. Fortunately, we can guess and verify their solutions by using a fixed-point technique combined with a regularization argument. We then show that ambiguity aversion to misspecifications should not be extremely large for the well-posedness of the optimization problems, at which point the existence of misspecified models is lost. Similar phenomena have been reported in previous studies of SDEs that are not bridges [43,44]. The proposed mathematical analysis scheme is technically difficult because it involves analyzing a nonlinear ODE with an exploding coefficient and a possibly unbounded solution. In this work, we focus on the CIR bridge (1); however, the theory developed here can be applied to other diffusion bridges with modifications.

A remarkable property of our formulation is that the worst-case misspecified bridges obtained from the corresponding pair of optimization problems explained above also have the form of (1) with different coefficient shapes of $r_t$ in their drift term. This implies that the proposed method does not breakdown the useful properties (well-posedness, numerical computability, and analytical tractability) of CIR bridges. A practical advantage of the proposed framework compared to the existing uncertainty quantification methods, such as the neural method [45] and probability density evolution method [46], is its high analytical tractability, as suggested above, where uncertainties can be quantified only through a system of ODEs, which can be computed easily; it is therefore free from any (possibly time-consuming statistical computations).

Finally, we apply misspecified CIR bridges to the estimation of daily fish counts from 10-min high-resolution fish count data for the well-studied model fish species *Plecoglossus altivelis* (*P. altivelis*) in Japan, which is a major inland fishery resource in the country, and its population dynamics have therefore received attention [e.g., 47,48,49,50,51]. We will discuss how the model performance depends on the degree of misspecification.

We close this subsection with a limitation of this study; the noise intensity is assumed not to be uncertain, while in applications, this parameter intensity may also become a source of uncertainties. In such a case, the Girsanov transformation will have a more complex form. Alternatively, there may be cases where the Girsanov transformation is simply inapplicable. For instance, if the Radon–Nikodym derivative between the benchmark model and the distorted model cannot be defined, a completely different optimization approach would likely be required. We intend to address these advanced issues in future research.

### 1.3 Structure of this paper

**Section 2** presents the diffusion bridge considered in this study and its basic properties. **Section 3**

formulates and theoretically examines a pair of optimization problems to estimate the cumulative average value of the solution to the diffusion bridge, which is a problem that arises in resource management. **Section 4** applies the diffusion bridge with possible misspecifications to a high-resolution fish count dataset in fisheries science. **Section 5** presents the summary and perspectives of this study. The **Appendix** presents technical proofs.

Overview of this study

① SDE (1) with uncertainties is given.
→ This is a "bridge-type" SDE.
② The total fish count $F$ is formulated based on the SDE.
→ This is an expectation functional.
③ Evaluate the upper- and lower- bounds of $F$.
→ A pair of optimization problems is specified.
④ Solve the optimality equations that give these bounds.
→ It reduces to solving some ODEs.
⑤ We also give a demonstrative application example.
→ The computation is due to solving these ODEs.

**Figure 1.** Overview diagram of the proposed framework.

## 2. Diffusion bridge and its misspecification

We formulate the target SDE and its uncertain counterparts. They are bridge-type SDEs aving certain unbounded coefficients.

### 2.1 Formulation

Throughout this paper, we work on a complete probability space $(\Omega, \mathcal{F}, \mathbb{P})$, where $\Omega$ is the set of all possible realizations, $\mathcal{F}$ is a collection of realizations in $\Omega$, and $\mathbb{P}$ is a probability measure. The expectation on $\mathbb{P}$ is denoted as $\mathbb{E}_{\mathbb{P}}$, and the same presentation rule is applied to the other probability measures.

The CIR bridge considered in this study is the following nonlinear SDE (reprint of (1)):

$$\underbrace{\mathrm{d}X_t}_{\text{Increment}} = \underbrace{\left( \underbrace{a_t}_{\text{Source}} - \underbrace{\frac{r_t}{1-t} X_t}_{\text{Reversion}} \right) \mathrm{d}t}_{\text{Drift}} + \underbrace{\underbrace{\sigma_t \sqrt{\frac{r_t}{1-t} X_t}}_{\text{Noise intensity}} \underbrace{\mathrm{d}B_t}_{\text{Noise increment}}}_{\text{Diffusion}}, \quad 0 < t < 1 \tag{2}$$

under the homogeneous initial and terminal conditions $X_0 = X_1 = 1$. Here, $B = (B_t)_{t \geq 0}$ is a 1-D standard Brownian motion, the diffusion term is defined in Itô's sense, and the coefficients $a, r, \sigma$ are time - dependent, continuous, positive, and bounded in the time interval $[0,1]$. Despite having exploding coefficients, the CIR bridge (2) admits a unique pathwise solution that is continuous and nonnegative (Proposition 1 in Yoshioka [52]) and hence is a well-posed mathematical model. This is due to the specific speed of the explosion in the drift and diffusion terms, where the explosion of the solution itself is prevented, and the terminal condition is continuously reached. Moreover, moments of the solution can be found in closed forms and are obtained explicitly in simple cases, such as the case where all coefficients $a, r, \sigma$ are positive constants (Appendix in Yoshioka [34]).

The CIR bridge (2) and its mean-field counterpart have been applied to modeling diurnal fish migration, where the solution $X$ represents the unit-time fish count observed at a fixed point in the river, and the normalized times 0 and 1 are identified as sunrise and sunset, respectively [52]. In these studies, it was assumed that the bridge can be obtained from observation data without any errors, which is an idealized condition because real observation projects are usually constrained by monetary and labor costs; hence, the quality and quantity of data collected are limited [53]. This motivated us to study the misspecified version of (2). Note that studies on the misspecified versions of diffusion bridges are limited thus far.

***Remark 1.*** Diffusion bridges with time-dependent diffusion coefficients also appear in machine learning, where the explosion of variance is efficiently managed [54,55]. In contrast, our bridge has an exploding diffusion coefficient but bounded variance, where the exploding drift and diffusion coefficient come from the consideration of a biological clock such that migration occurs during sunrise and sunset. Therefore, an unbounded diffusion coefficient does not necessarily lead to unbounded variance. In our case, the unbounded reversion prevents the explosion of variance of $X$.

### 2.2 Modeling uncertainties

Hansen and Sargent [25] assumed that model uncertainty is represented by an added drift term in an SDE. Their framework adapted for our formulation is explained in the following section. The explanation below is based on Baltas and Yannacopoulos [33] with an adaptation to our setting.

First, consider a real scalar process $u = (u_t)_{0<t<1}$ that is adapted to natural filtration $\mathcal{F}$ (such a process is called an adapted process hereafter) generated by Brownian motion $B$ such that the following Novikov condition is satisfied (see also Section 2 in Hess [56]):

$$\mathbb{E}_{\mathbb{P}}\left[\exp\left(\frac{1}{2}\int_0^1 u_s^2 \mathrm{d}s\right)\right] < +\infty . \tag{3}$$

This condition means that the process $u$ is sufficiently regular in the sense that its quadratic integration is small. Under the Novikov condition (3), given one $u$, we can define a probability measure $\mathbb{Q}(u)$ that is equivalent to $\mathbb{P}$, where the Radon–Nikodym derivative between them is given as follows:

$$\left.\frac{\mathrm{d}\mathbb{Q}(u)}{\mathrm{d}\mathbb{P}}\right|_t = \exp\left(\int_0^t u_s \mathrm{d}B_s - \frac{1}{2}\int_0^t u_s^2 \mathrm{d}s\right), \quad 0 < t < 1, \tag{4}$$

which is a martingale. We find that the following process $W = (W_t)_{0\le t\le 1}$ is a 1-D standard Brownian motion on the new probability measure $\mathbb{Q}(u)$ (i.e., Girsanov transformation):

$$W_t = B_t - \int_0^t u_s \mathrm{d}s . \tag{5}$$

With (5) in mind, our CIR bridge on $\mathbb{Q}(u)$ is given by

$$\mathrm{d}X_t = \left( a_t - \frac{r_t}{1-t} X_t \underbrace{+\sigma_t u_t \sqrt{\frac{r_t}{1-t} X_t}}_{\text{Added drift coefficient}} \right) \mathrm{d}t + \sigma_t \sqrt{\frac{r_t}{1-t} X_t} \mathrm{d}W_t, \quad 0 < t < 1, \tag{6}$$

where the added drift coefficient is proportional to $u_t$. The initial and terminal conditions remain unchanged. Finally, the relative entropy $R(u)$ between the two probability measures $\mathbb{P}$ and $\mathbb{Q}(u)$ in our formulation is obtained as follows:

$$R(u) = \mathbb{E}_{\mathbb{P}}\left[ \left.\frac{\mathrm{d}\mathbb{Q}(u)}{\mathrm{d}\mathbb{P}}\right|_1 \ln\left( \left.\frac{\mathrm{d}\mathbb{Q}(u)}{\mathrm{d}\mathbb{P}}\right|_1 \right) \right] = \mathbb{E}_{\mathbb{Q}(u)}\left[ \ln\left( \left.\frac{\mathrm{d}\mathbb{Q}(u)}{\mathrm{d}\mathbb{P}}\right|_1 \right) \right] = \frac{1}{2}\mathbb{E}_{\mathbb{Q}(u)}\left[ \int_0^1 u_s^2 \mathrm{d}s \right], \quad 0 < t < 1. \tag{7}$$

Note that the expectation is evaluated on $\mathbb{Q}(u)$ on the right-most side. This is a quadratic functional of $u$ that often appears in control problems, which is the reason that the framework of Hansen and Sargent [25] harmonizes well with stochastic control formalism. We have $R(u) = 0$ when $u \equiv 0$ (i.e., when there is no misspecification).

We close this section with a discussion about the meaning and implications of the model of uncertainties presented above. In the present framework, the probability measure $\mathbb{P}$ represents a benchmark model, e.g., a model identified from data. The other probability measure $\mathbb{Q}(u)$ given a suitable $u$ represents a distorted model where the difference between the two probability measures is defined through the Radon–Nikodym derivative (4) and is quantified by the relative entropy (7). In the context of SDEs, the distorted model (6) has an added drift coefficient whose influence increases as $u$ increases in the sense that the relative entropy (7) increases. Because we assume the existence of a Radon–Nikodym derivative between equivalent probability measures $\mathbb{P}$ and $\mathbb{Q}(u)$, their difference is constrained at least quantitatively. For example, the satisfaction or violation of the Feller condition (i.e., whether the solution to an SDE of the CIR type hits 0 or not) should be preserved between equivalent probability measures [57]. In this view, the functional form of $u$, which possibly depends implicitly on $X$, is not completely arbitrary. In some tractable SDEs, the condition of the existence of equivalence between benchmark and distorted models can be theoretically characterized [58]. We will also obtain such a condition in the next section.

Note that the distorted model (6) may be fundamentally different from the benchmark model (2) owing to the added drift coefficient, which is apparently nonlinear with respect to $X_t$; hence, the well-posedness of (6) on $\mathbb{Q}(u)$ is not trivial at this stage. In the next section, we show that, in our case, this issue can be effectively resolved by measuring model uncertainties through relative entropy.

***Remark 2.*** Misspecifications in the context of Hansen and Sargent [25] that we employ in this paper do not change the Feller condition between benchmark and misspecified models. Indeed, this condition depends on the diffusion term and the source $a$ (e.g., Yoshioka [34]), but the misspecified term we will derive in the next section turns out to affect only the reversion term in the drift coefficient.

## 3. Upper and lower bounding models

We formulate the target optimization problems and mathematically study them. We see that these problems reduce to some pair of ODEs.

### 3.1 Overview

We consider the worst-case under- and overestimation problems of the following integral quantity $F$:

$$F = \mathbb{E}_{\mathbb{P}}\left[\int_0^1 X_s \mathrm{d}s\right], \tag{8}$$

which, in the context of the fish count problem, is the expectation of the daily total fish count because $X$ represents the unit-time fish count and the time interval $[0,1]$ is from sunrise to sunset, during which fish migrate; they are assumed to be sleeping during nighttime. If the model used to evaluate $F$ is distorted, then the expectation must be redefined using $\mathbb{E}_{\mathbb{Q}(u)}$ with a suitable $u$. Both under- and overestimations of $F$ are theoretically possible. Therefore, we consider a pair of optimization problems to estimate them under the constraint that the size of the uncertainties is constrained in terms of the relative entropy.

### 3.2 Lower-bounding case

We formulate a lower-bounding problem of the expectation $F$ in (8) considering possible misspecifications. For an adapted process $u$ satisfying the Novikov condition (3), we set the following minimization problem: find

$$G(t,x) = \inf_u \left( \underbrace{\mathbb{E}_{\mathbb{Q}(u)}\left[\int_t^1 X_s \mathrm{d}s \middle| X_t = x\right]}_{\text{Target quantity}} + \underbrace{\frac{1}{2\psi}\mathbb{E}_{\mathbb{Q}(u)}\left[\int_t^1 u_s^2 \mathrm{d}s \middle| X_t = x\right]}_{\text{Relative entropy}} \right), \quad 0 \le t \le 1, \quad x \ge 0 \tag{9}$$

subject to the dynamics (6). Here, we assume that this SDE (6) admits a pathwise unique solution that attains the terminal condition $X_1 = 0$ continuously in time (i.e., $\lim_{t \to 1-0} X_t = 0$ with probability 1). This is a minimization problem of an objective function containing the target expectation and the relative entropy, whose contribution is modulated by the ambiguity-aversion parameter $\psi > 0$ that plays the role of a Lagrangian multiplier. A larger value of the (time-restricted) relative entropy (9) is allowed to become larger as $\psi$ increases, and the limit $\psi \to 0$ reduces to the ambiguity-neutral case where the decision-maker of the optimization problem, e.g., an observer of the fish count, ignores model uncertainties.

In view of the definition (8) of $F$, the quantity of interest in application is the special case $G(0,0)$, but considering the problem of the form (9) allows us to apply a dynamic programming argument with which we can reduce the problem (9) to finding an appropriate solution to an optimality equation that is solvable analytically. Finally, a minimizing $u$ of (9) is denoted by $u^*$. We have $G(0,0) \le F$ by substituting the (admissible) null control $u \equiv 0$ into the quantity to be minimized in (9).

### 3.3 Upper-bounding case

The upper-bounding case is formulated symmetrically by considering the corresponding lower-bounding case as follows:

$$H(t,x) = \sup_u \left( \underbrace{\mathbb{E}_{\mathbb{Q}(u)}\left[ \int_t^1 X_s \mathrm{d}s \middle| X_t = x \right]}_{\text{Target quantity}} - \frac{1}{2\psi} \underbrace{\mathbb{E}_{\mathbb{Q}(u)}\left[ \int_t^1 u_s^2 \mathrm{d}s \middle| X_t = x \right]}_{\text{Relative entropy}} \right), \quad 0 \le t \le 1, \quad x \ge 0 \tag{10}$$

subject to the dynamics (6). Here, we assume that this SDE (6) admits a pathwise unique solution that attains the terminal condition $X_1 = 0$ continuously in time. We have $H(0,0) \ge F$ by substituting the (admissible) null control $u \equiv 0$ into the quantity to be minimized in (10).

***Remark 3.*** The optimization problems presented in this subsection are unconstrained optimization problems that can be seen as relaxed versions of the constrained ones, where the ambiguity-aversion parameter $\psi$ serves as a Lagrangian multiplier (e.g., Melo et al. [59]). For the lower-bounding case, we can rewrite (9) as follows:

$$G(t,x) = \inf_u \left( \mathbb{E}_{\mathbb{Q}(u)}\left[ \int_t^1 X_s \mathrm{d}s \middle| X_t = x \right]; R(u) = \frac{1}{2} \mathbb{E}_{\mathbb{Q}(u)}\left[ \int_t^1 u_s^2 \mathrm{d}s \middle| X_t = x \right] \le \kappa \right), \quad 0 \le t \le 1, \quad x \ge 0, \tag{11}$$

where $\kappa > 0$ is a parameter that gives the size of the model uncertainties. We can formally consider $\psi$ as a function of $\kappa$, that is, $\psi = \psi(\kappa)$. From this perspective, by assuming that $\psi(\kappa)$ is invertible, the unconstrained problem (9) becomes equivalent to (11). The same reasoning applies to the upper-bounding case. The optimization problems (9) and (10) therefore give the worst-case lower and upper bounds of $F$ under a given uncertainty size $\kappa$. We exploit these relationships in the application study described in **Section 4**.

***Remark 4.*** We have the following inequality, which is the reason for naming the "lower-bounding" and "upper-bounding" cases:

$$G(0,0) \le F \le H(0,0), \quad \psi > 0. \tag{12}$$

### 3.4 Mathematical analysis

In this subsection, we study the upper-bounding case (10) exclusively because the lower-bounding case (9) is much easier in theory (**Remark 4**). More specifically, solving the upper-bounding case involves dealing with a nonlinear ODE with an exploding coefficient, whose solution possibly blows up in finite time. In contrast, this difficulty is not encountered in the lower-bounding case. Moreover, we focus on the case where the reversion parameter is a constant $r_t \equiv r > 0$ because this is the situation that we consider in **Section 4**.

The following proposition (**Proposition 1**) is our main theoretical result, which states that the optimization problem (10) can be solved if $\psi > 0$ is sufficiently small. We set $\bar{\sigma} = \max_{0 \le t \le 1} |\sigma_t| > 0$.

***Proposition 1***

*Assume that*

$$\bar{\sigma}^2 \psi < \frac{r}{2}. \tag{13}$$

*Then, it follows that*

$$H(t,x) = A_t x + B_t, \quad 0 \le t \le 1, \quad x \ge 0, \tag{14}$$

*where*

$$B_t = \int_t^1 a_s A_s \mathrm{d}s \tag{15}$$

*and* $A = (A_t)_{0 \le t \le 1}$ *is the unique continuously differentiable solution to the following time-backward ODE:*

$$-\frac{\mathrm{d}A_t}{\mathrm{d}t} = \frac{r}{1-t}\left(-A_t + \frac{\sigma_t^2 \psi}{2} A_t^2\right) + 1, \quad 0 \le t < 1, \quad A_1 = 0. \tag{16}$$

*Moreover, it follows that the maximizer in* (10) *is given by*

$$u_t^* = A_t \sigma_t \psi \sqrt{\frac{r}{1-t} X_t}, \quad 0 \le t < 1 \tag{17}$$

*and the controlled process* $X$ *is governed by*

$$\mathrm{d}X_t = \left(a_t - \frac{r}{1-t}\left(1 - \sigma_t^2 \psi A_t\right) X_t\right)\mathrm{d}t + \sigma_t \sqrt{\frac{r}{1-t} X_t}\,\mathrm{d}W_t, \quad 0 < t < 1 \tag{18}$$

*with* $X_0 = X_1 = 0$.

The main implication of **Proposition 1** is that we can find $u_t^*$ in a closed form and that the controlled process is again a CIR bridge of the form (2). This implies that, under the assumed model uncertainties, the form of the target system dynamics does not formally change. From a practical standpoint, this means that analytical and computational methodologies to address bridges (2) carry over to (18), thereby facilitating the application of a misspecified model. The most technical part of **Proposition 1** is the ODE (16). Indeed, this ODE may have a solution only locally in time, whereas for the well-posedness of the optimization problem, we need the solution to be globally defined in the time interval $[0,1]$. The assumption (13), which means that the ambiguity-aversion parameter needs to be sufficiently small (i.e., model uncertainties should not be assumed to be too large), guarantees this condition. The assumption (13) is a sufficient condition, but not a necessary and sufficient condition, implying that it is not sharp. Nevertheless, it is suggested that the ambiguity-aversion parameter $\psi$ must be sufficiently small and is computationally supported in our application study in **Section 4**.

***Remark 5.*** The lower-bounding case can be solved in the same way as in **Proposition 1** with the following

modifications:

$$-\frac{\mathrm{d}A_t}{\mathrm{d}t}=\frac{r}{1-t}\left(-A_t-\frac{\sigma_t^2\psi}{2}A_t^2\right)+1,\quad 0\le t<1,\quad A_1=0 \tag{19}$$

and

$$u_t^*=-A_t\sigma\psi\sqrt{\frac{r}{1-t}X_t},\quad 0\le t<1, \tag{20}$$

and the controlled process satisfies the SDE

$$\mathrm{d}X_t=\left(a_t-\frac{r}{1-t}\left(1+\sigma_t^2\psi A_t\right)X_t\right)\mathrm{d}t+\sigma_t\sqrt{\frac{r}{1-t}X_t}\mathrm{d}W_t,\quad 0<t<1. \tag{21}$$

## 4. Application

We present a demonstrative application example of the proposed approach. The data focused on here are 10-min upstream-migrating fish count.

### 4.1 Target data

We apply the proposed lower- and upper-bounding cases using the CIR bridge under the misspecifications studied in **Section 3** to 10-min, AI-based high-resolution fish count data[1] of juvenile *P. altivelis* in spring collected at the Nagara River Barrage near the mouth of the Nagara River (Class-A River in the Tokai Region) in Japan. The Nagara River is a major river to which *P. altivelis* migrates every year. We have the latest dataset of the 10-min fish count from approximately the end of February to the end of June in 2023, 2024, and 2025 courtesy of the Japan Water Agency, Ibi River, and Nagara River General Management Office. Counting the migrating population of juvenile *P. altivelis* has been considered an important task for assessing its stock in inland waters in Japan. However, high-resolution data such as those collected at the Nagara River are still not common and in most cases are counted manually with coarser quality and quantity, as noted by Yoshioka [53].

In a previous study, we identified the CIR bridge (2) by assuming constant parameters $a,r,\sigma>0$ based on a least-squares procedure using empirical and theoretical time-dependent averages and standard deviations [52]. The time interval between sunrise and sunset each day was normalized to the unit interval $(0,1)$ because juvenile *P. altivelis* migrates during the daytime, as discussed in Yoshioka [34], and the unit-time fish count $X$ was normalized using the daily fish count and data resolution (10 min). Detailed descriptions of the data and model identification methods have been presented in previous studies [34,52] and are not repeated here. One important finding was that the trajectories of the identified CIR bridge were nonnegative and intermittent, and the latter property may result in the loss of modeling accuracy

[1] Japan Water Agency, Ibi River and Nagara River General Management Office https://www.water.go.jp/chubu/nagara/15_sojou/chousahouhou.html (Last accessed on October 18, 2025)

if the fish are counted manually with a low temporal resolution [53]. The identified parameter values for the constant-parameter case were as follows (the right-most column of Table B.1 in Yoshioka [52]): $a = 0.03673$, $r = 0.7100$, and $\sigma = 0.7252$, which are used in this paper. A more complex model can also be employed when necessary, but here we focus on the simplest case for demonstration purposes.

Under the present setting, the only remaining parameter is the ambiguity-aversion parameter $\psi$, which appears in the optimization problems. Therefore, we investigate these problems for different values of $\psi$ and the performance of the benchmark and misspecified (i.e., lower- and upper-bounding) models.

### 4.2 Results and discussion

We first discuss the performance of the misspecified models using the expected total daily counts predicted for the lower- and upper-bounding cases. We assume that $G(0,0)$ and $H(0,0)$ exist by solving their corresponding optimality equations. For the lower-bounding case, we have

$$G(0,0) = \mathbb{E}_{\mathbb{Q}(u^*)}\left[\int_0^1 X_s \mathrm{d}s\right] + \frac{1}{2\psi}\mathbb{E}_{\mathbb{Q}(u^*)}\left[\int_0^1 \left(u_s^*\right)^2 \mathrm{d}s\right], \tag{22}$$

where the left-hand side of (22) is now available and the first term on the right-hand side of (22) is obtained by solving the ODE (see **Remark 3**):

$$\frac{\mathrm{d}}{\mathrm{d}t}\mathbb{E}_{\mathbb{Q}(u^*)}\left[X_s\right] = a - \frac{r}{1-t}\left(1+\sigma^2\psi A_t\right)\mathbb{E}_{\mathbb{Q}(u^*)}\left[X_s\right],\quad 0 < t < 1, \tag{23}$$

starting from the initial condition 0. We then obtain the corresponding relative entropy as follows:

$$\kappa(\psi) = \frac{1}{2}\mathbb{E}_{\mathbb{Q}(u^*)}\left[\int_0^1 \left(u_s^*\right)^2 \mathrm{d}s\right] = \psi\left(G(0,0) - \mathbb{E}_{\mathbb{Q}(u^*)}\left[\int_0^1 X_s \mathrm{d}s\right]\right), \tag{24}$$

considering **Remark 3** and the fact that the right-most side is now available. The ODE is linear and can be numerically discretized. We apply a naïve one-step method with a time increment $10^{-6}$, which implicitly discretizes the term multiplied by the exploding coefficient $\frac{r}{1-t}$. Essentially, the same theoretical and computational reasoning applies to the upper-bounding case, where the ODE to be solved is given by

$$\frac{\mathrm{d}}{\mathrm{d}t}\mathbb{E}_{\mathbb{Q}(u^*)}\left[X_s\right] = a - \frac{r}{1-t}\left(1-\sigma^2\psi A_t\right)\mathbb{E}_{\mathbb{Q}(u^*)}\left[X_s\right],\quad 0 < t < 1 \tag{25}$$

and the relative entropy by

$$\kappa(\psi) = \frac{1}{2}\mathbb{E}_{\mathbb{Q}(u^*)}\left[\int_0^1 \left(u_s^*\right)^2 \mathrm{d}s\right] = \psi\left(\mathbb{E}_{\mathbb{Q}(u^*)}\left[\int_0^1 X_s \mathrm{d}s\right] - H(0,0)\right). \tag{26}$$

An advantage of the proposed method in this application is that we do not need to simulate sample paths of the CIR bridges, which may be time consuming and introduce unnecessary statistical errors in the analysis.

The expectation of the total daily fish count for the benchmark case is $F = 0.01074$ with the prescribed parameter values (Appendix of Yoshioka [34]). Based on the methodology explained above, **Figure 2** shows the expected total daily counts as a function of the ambiguity-aversion parameter $\psi$ for the lower-bounding, benchmark, and upper-bounding cases. We find that the benchmark case is bounded by misspecified cases and that the overestimation case provides a rapidly increasing prediction of the total daily count. **Figure 3** shows the computed $A$ from the ODEs (16) and (19). In the present setting, we encountered the divergence toward $+\infty$ of numerical solutions to the ODE (16) near $\psi = 15.8$ for the upper-bounding case, while this problem did not occur for the lower-bounding case. This computationally supports the assumption (13) in **Proposition 1** that the ambiguity-aversion parameter $\psi$ must be chosen to be sufficiently small so that the upper-bounding problem is well-posed. The upper bound of $\psi$ obtained from (13) is $\frac{r}{2\sigma^2} = 0.68 < 15.8$ and is therefore conservative because it is a sufficient condition but not a necessary and sufficient condition. Below, the values of $\psi$ are chosen so that this condition is satisfied.

**Figure 4** shows the relative entropy size $\kappa = \kappa(\psi)$ as a function of the ambiguity-aversion parameter $\psi$ for misspecified cases. The increase in $\kappa$ for the upper-bounding case is sharper than that of the underestimation case, reflecting the explosion of the worst-case total daily fish count in the former case reported in **Figure 2**. For both cases, the curves of $\kappa = \kappa(\psi)$ are increasing, convex, and vanishing at $\psi = 0$, supporting its invertibility.

**Figures 5 and 6** show the expected total daily counts as functions of the relative entropy size for the lower- and upper-bounding cases, respectively. Thus, these figures visualize the expected under- and overestimations of the total daily fish counts as functions of the size of uncertainties, which correspond to solving constrained optimization problems that bound the benchmark case (**Remark 3**). In both cases, the misspecified models yield more extreme results as the relative entropy increases. A practical interest here is the possible size of the uncertainties, namely, the values of $\kappa$, in applications. Yoshioka [53] evaluated the juvenile migration of *P. altivelis* and found that conventional manual observations would under- and overestimate the total daily fish count with relative errors of the same order as the true values. These evaluation results imply that the quantity $\mathbb{E}_{\mathbb{Q}(u^*)}\left[\int_0^1 X_s \mathrm{d}s\right] / F$ in the order of $10^{-1}$ to $10^1$ is possible. In the present case, **Figures 5 and 6** suggest that the possible maximum values of $\kappa$ for the lower- and upper-bounding cases are approximately 0.2 to 0.5 and 1 to 3 (**Tables 1 and 2**), respectively, which are one order different from each other. We consider that these order estimates would apply to the migration of *P. altivelis* at other study sites if datasets become available, while their applicability to other fish species is unclear because the migration patterns are highly dependent on the species. Nevertheless, as demonstrated in this subsection, we can explicitly link the relative entropy and performance of the misspecified models owing to the tractability of the CIR bridge and the pair of optimization problems.

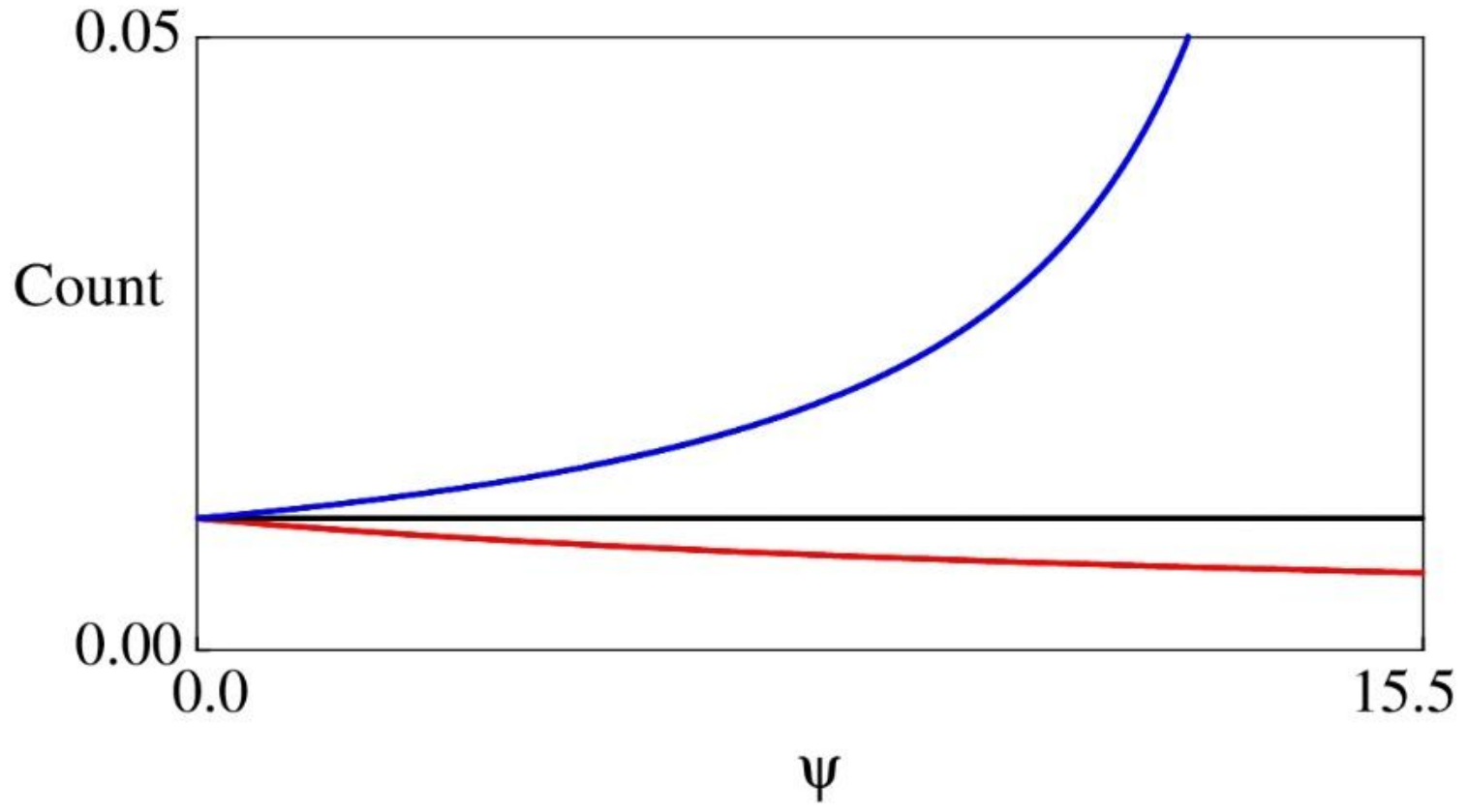

**Figure 2.** Expected total daily count ("Count") as a function of the ambiguity-aversion parameter $\psi$ : lower-bounding case (red), benchmark case (black), and upper-bounding case (blue).

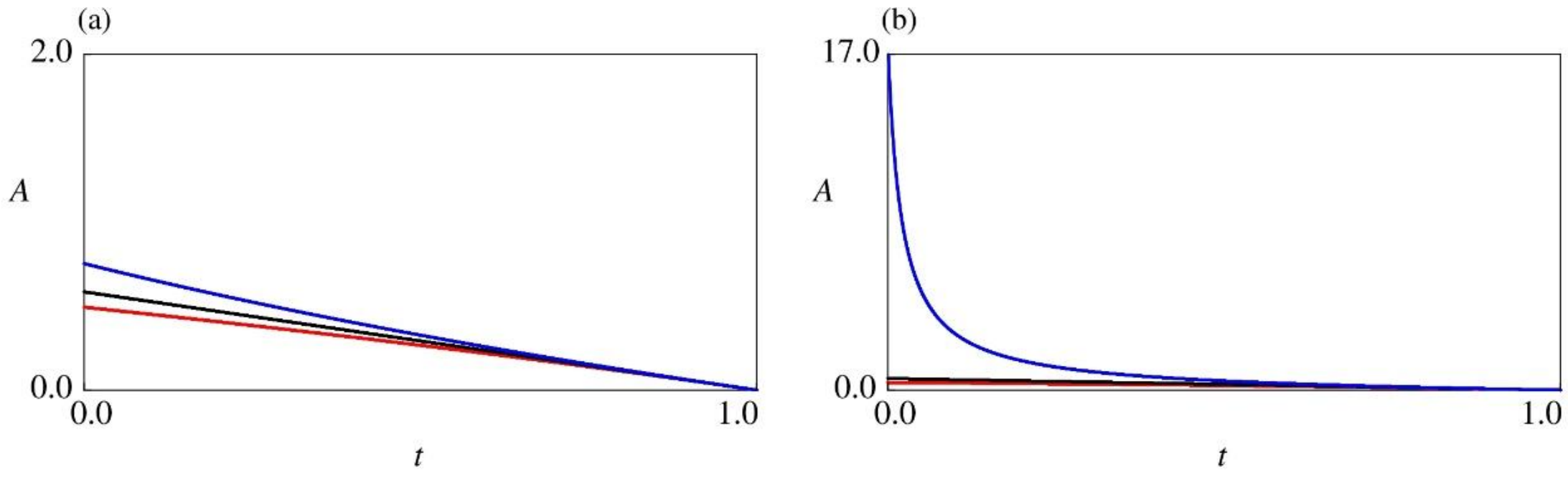

**Figure 3.** Computed $A$ for (a) $\psi = 5.0$ and (b) $\psi = 15.5$ : lower-bounding case (red), benchmark case (black), and upper-bounding case (blue). Note the difference in vertical scales between (a) and (b).

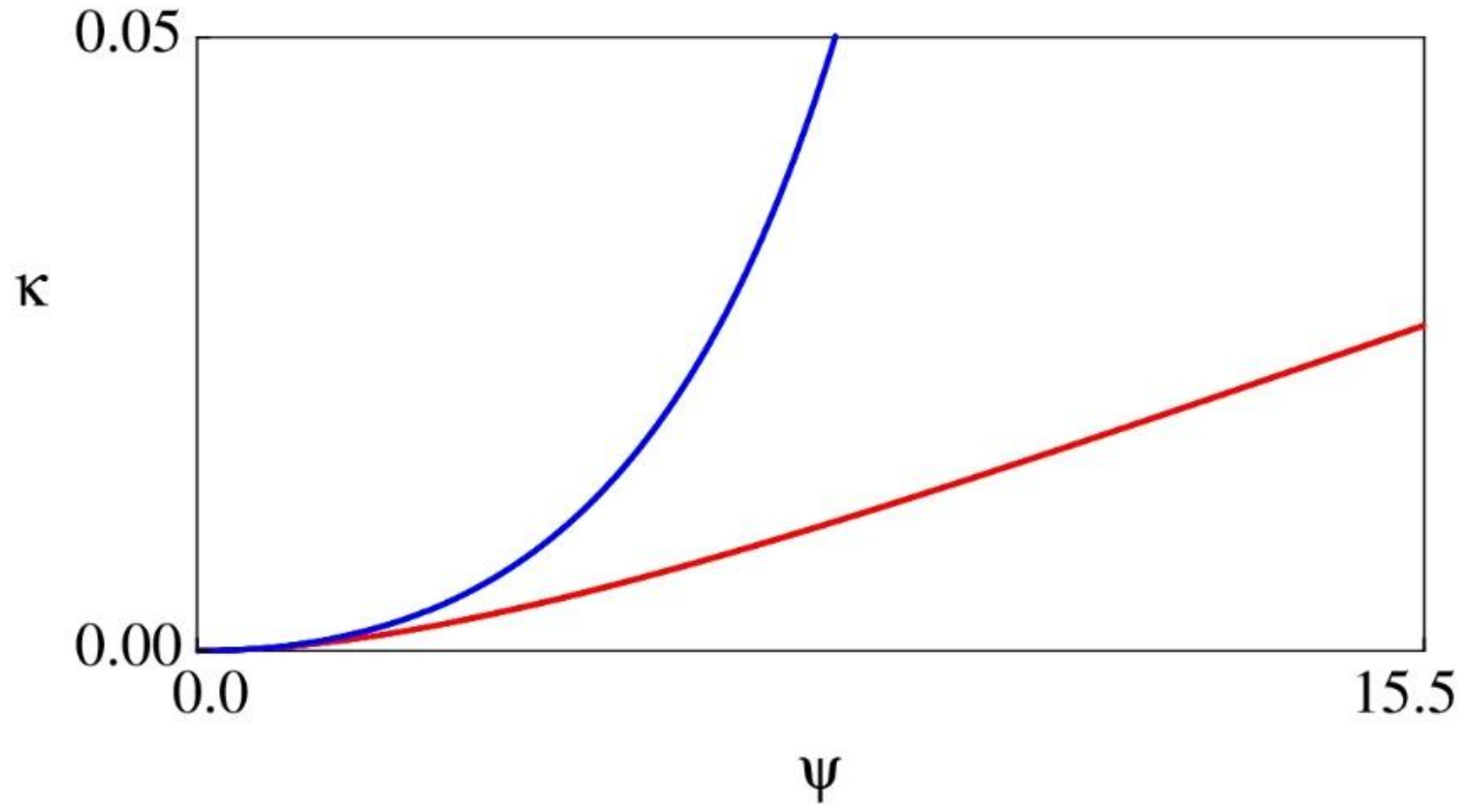

**Figure 4.** Relative entropy size $\kappa$ , i.e., size of uncertainties, as a function of ambiguity-aversion parameter $\psi$ : lower-bounding case (red) and upper-bounding case (blue).

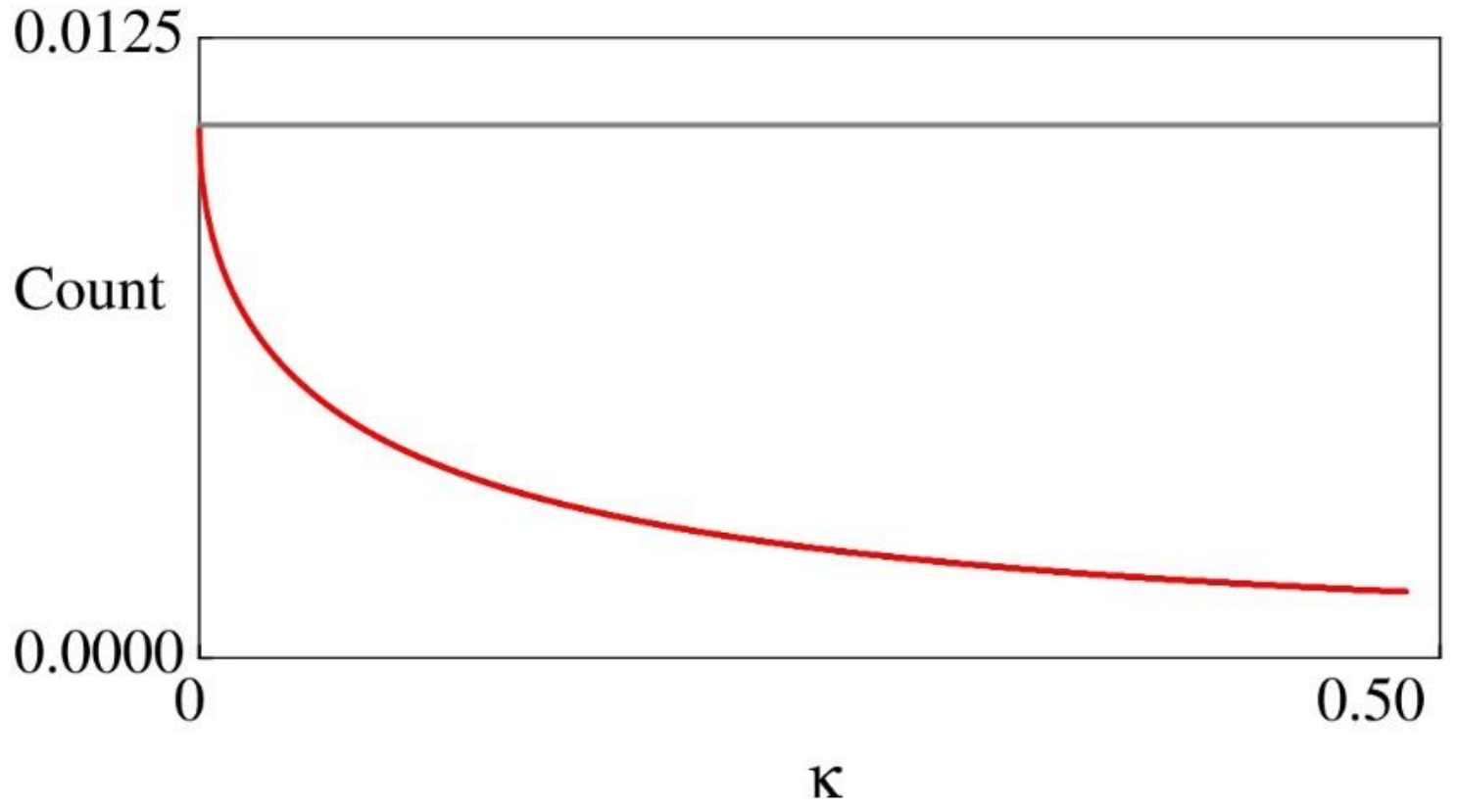

**Figure 5.** Total count in the lower-bounding case as a function of relative entropy size $\kappa$. Gray line for the benchmark case.

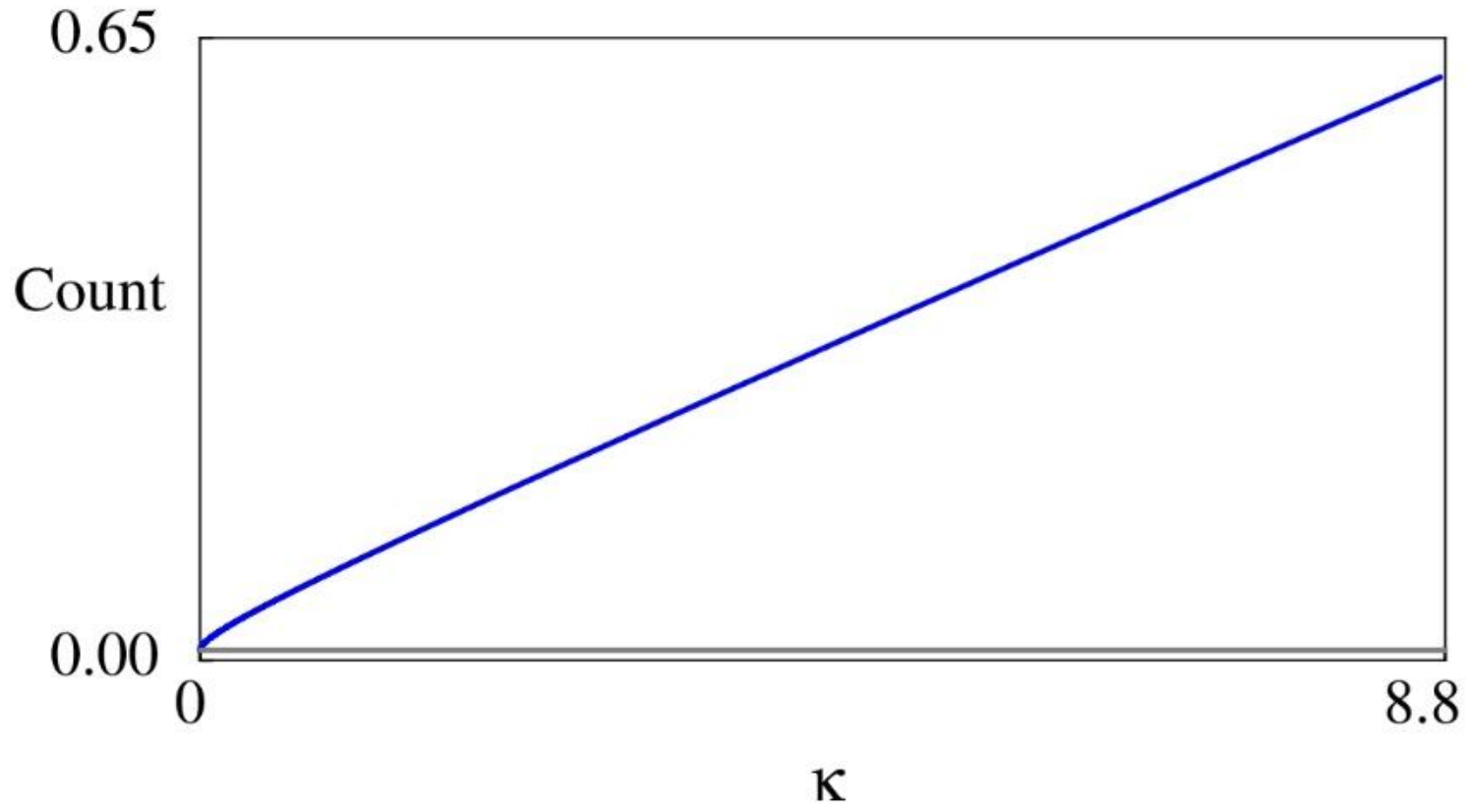

**Figure 6.** Total count in the upper-bounding case as a function of relative entropy size $\kappa$. Gray line for the benchmark case.

**Table 1.** Values of the count ratio $\mathbb{E}_{\mathbb{Q}(u^*)}\left[\int_0^1 X_s \mathrm{d}s\right]/F$ for different values of $\psi$ and $\kappa$ in the lower-bounding case.

| $\psi$ | $\kappa$ | Count ratio |
|---|---|---|
| 5 | 4.91.E-03 | 7.98.E-01 |
| 10 | 1.45.E-02 | 6.76.E-01 |
| 50 | 9.64.E-02 | 3.61.E-01 |
| 100 | 1.73.E-01 | 2.60.E-01 |
| 400 | 4.51.E-01 | 1.32.E-01 |

**Table 2.** Values of the count ratio $\mathbb{E}_{\mathbb{Q}(u^*)}\left[\int_0^1 X_s \mathrm{d}s\right]/F$ for different values of $\psi$ and $\kappa$ in the upper-bounding case.

| $\psi$ | $\kappa$ | Count ratio |
|---|---|---|
| 5 | 1.22.E-02 | 1.40.E+00 |
| 10 | 1.13.E-01 | 2.57.E+00 |
| 13 | 4.84.E-01 | 5.48.E+00 |
| 14 | 9.62.E-01 | 8.75.E+00 |
| 15 | 2.83.E+00 | 2.06.E+01 |

## 5. Conclusion

We studied the lower- and upper-bounding misspecifications of the integral of the solutions to a CIR bridge with exploding drift and diffusion coefficients. The worst-case misspecifications introduced using the measure transformation and relative entropy were preserved with the well-posed and analytical tractability of the original CIR bridge while handling a wide range of uncertainty sizes. A fishery case study based on the proposed framework for misspecifications suggested how the ambiguity-aversion parameter can be prescribed in applications. Consequently, this study contributes to the formulation, analysis, and applications of diffusion bridges. To the best of the author's knowledge, minute-resolution datasets of fish migration are still rare. As more high-resolution fish migration datasets become available, it is expected that the modeling of fish migration itself, as well as the assessment of its uncertainties, will be advanced.

We consider that the proposed framework also applies to diffusion bridges studied in other research areas [60,61], although their misspecified versions need to be investigated computationally if the corresponding optimality equations do not have analytical solutions and/or if the underlying stochastic processes are more complex than those studied in this paper. In this view, the contribution of this study will have a spillover effect on research fields other than fishery science. A theoretical issue that was not addressed in this study is the study of other uncertainty forms, such as those modeled using generalized divergences [62] that generalize relative entropy. A practical drawback in this case is the increase in the number of parameters for modeling uncertainties, although problems with generalized divergences would be a theoretically interesting research topic. The application of the proposed framework to multidimensional and non-Markovian SDEs and diffusion bridges is also an interesting topic.

## Appendix

Before proceeding to **Proof of Proposition 1**, we prepare **Lemma 1**, which plays a role in the proof. The right-hand side of inequality (27) actually has redundancy but is presented in that form for tractability in the **Proof of Lemma 1**. We note that $\frac{r}{2}$ is the smallest on the right-hand side of (27).

***Lemma 1***

*Assume that (reprint of (13))*

$$\bar{\sigma}^2\psi < \min\left\{ r+1, \frac{(r+1)^2}{2r}, \frac{2(r+1)^2}{r}, \frac{(r+1)(r+4)}{2(r+2)}, \frac{r}{2} \right\}. \quad (27)$$

*Then, it is possible to choose a constant* $\bar{Y} > 0$ *such that the following four conditions are satisfied:*

$$\frac{r}{(r+1)(r+2)}\bar{\sigma}^2\psi(\bar{Y}+1) < 1, \quad (28)$$

$$\frac{r}{2(r+1)(r+2)}\bar{\sigma}^2\psi(\bar{Y}+1)^2 < \bar{Y}, \quad (29)$$

$$\frac{1}{r+1}\bar{\sigma}^2\psi(\bar{Y}+1) < 1, \quad (30)$$

*and*

$$\frac{r}{2(r+1)^2}\bar{\sigma}^2\psi(\bar{Y}+1)^2 < 1. \quad (31)$$

***Proof of Lemma 1***

For later use, set

$$Q = \frac{r}{(r+1)(r+2)}\bar{\sigma}^2\psi > 0. \quad (32)$$

Then, we have

$$Q < \frac{1}{2} \quad (33)$$

due to (27):

$$Q = \frac{r}{(r+1)(r+2)}\bar{\sigma}^2\psi < \frac{r}{(r+1)(r+2)}\frac{(r+1)^2}{2r} = \frac{r+1}{2(r+2)} < \frac{1}{2}. \quad (34)$$

By (32), the first condition (28) is rewritten as

$$Q(\bar{Y}+1) < 1 \leftrightarrow \bar{Y} < \frac{1}{Q} - 1. \quad (35)$$

The right-hand side of (35) is positive because of (33): similarly, the third condition (30) is rewritten as

$$\frac{r+2}{r}Q(\bar{Y}+1) < 1 \leftrightarrow \bar{Y} < \frac{r}{r+2}\frac{1}{Q} - 1 \left( = \frac{r+1}{\bar{\sigma}^2\psi} - 1 \right). \quad (36)$$

The right-hand side of the second inequality in (36) is positive because (27). The second condition (29) can be rewritten as

$$\frac{Q}{2}\left(\bar{Y}+1\right)^2<\bar{Y}\,, \tag{37}$$

which is a quadratic inequality of $\bar{Y}$ solvable as follows (recall (33)):

$$0<\frac{1}{Q}-1-\frac{\sqrt{1-2Q}}{Q}<\bar{Y}<\frac{1}{Q}-1+\frac{\sqrt{1-2Q}}{Q}\,. \tag{38}$$

Note that the left-most inequality in (38) is true because

$$\left(\frac{1}{Q}-1\right)^2-\left(\frac{\sqrt{1-2Q}}{Q}\right)^2=\frac{1}{Q^2}\left(\left(1-Q\right)^2-\left(1-2Q\right)\right)=1>0\,. \tag{39}$$

The fourth condition (31) is rewritten as

$$\frac{r+2}{2\left(r+1\right)}Q\left(\bar{Y}+1\right)^2<1\,, \tag{40}$$

which can be solved for $\bar{Y}$ as

$$\bar{Y}<\sqrt{\frac{2\left(r+1\right)}{r+2}\frac{1}{Q}}-1\;\left(=\sqrt{\frac{2\left(r+1\right)^2}{r}\frac{1}{\bar{\sigma}\psi^2}}-1\right). \tag{41}$$

The right side of (41) is positive by (27):

We can choose $\bar{Y}>0$ such that both the first condition (28) and third condition (29) are satisfied because of (35) and (38):

$$\frac{1}{Q}-1-\frac{\sqrt{1-2Q}}{Q}<\bar{Y}<\frac{1}{Q}-1\,. \tag{42}$$

Due to (36), the second condition is further satisfied if

$$\frac{1}{Q}-1-\frac{\sqrt{1-2Q}}{Q}<\bar{Y}<\frac{r}{r+2}\frac{1}{Q}-1\,. \tag{43}$$

This inequality makes sense because of

$$\frac{r}{r+2}\frac{1}{Q}-\left(\frac{1}{Q}-\frac{\sqrt{1-2Q}}{Q}\right)=\frac{\sqrt{1-2Q}}{Q}-\frac{2}{r+2}\frac{1}{Q}=\frac{1}{Q}\left(\sqrt{1-2Q}-\frac{2}{r+2}\right) \tag{44}$$

and

$$\left(\sqrt{1-2Q}\right)^2-\left(\frac{2}{r+2}\right)^2 =1-2Q-\frac{4}{\left(r+2\right)^2}$$
$$=\frac{\left(r+2\right)^2-4}{\left(r+2\right)^2}-\frac{2r}{\left(r+1\right)\left(r+2\right)}\bar{\sigma}^2\psi$$
$$=\frac{r\left(r+4\right)}{\left(r+2\right)^2}-\frac{2r}{\left(r+1\right)\left(r+2\right)}\bar{\sigma}^2\psi \qquad , \tag{45}$$
$$=\frac{2r}{\left(r+1\right)\left(r+2\right)}\left(\frac{\left(r+1\right)\left(r+4\right)}{2\left(r+2\right)}-\bar{\sigma}^2\psi\right)$$
$$>0$$

where we used (27). The remaining task is to demonstrate that we can choose $\bar{Y}>0$ such that (41) and (43) are satisfied. To show this, we prove the following inequality:

$$\frac{1}{Q}-1-\frac{\sqrt{1-2Q}}{Q}<\sqrt{\frac{2\left(r+1\right)}{r+2}\frac{1}{Q}}-1<\frac{r}{r+2}\frac{1}{Q}-1. \tag{46}$$

The right inequality of (46) follows from

$$\left(\frac{r}{r+2}\frac{1}{Q}\right)^2-\left(\sqrt{\frac{2\left(r+1\right)}{r+2}\frac{1}{Q}}\right)^2=\frac{r^2}{\left(r+2\right)^2}\frac{1}{Q^2}-\frac{2\left(r+1\right)}{\left(r+2\right)Q}$$
$$=\frac{r^2}{\left(r+2\right)^2}\left(\frac{1}{Q^2}-\frac{2\left(r+1\right)\left(r+2\right)}{r^2}\frac{1}{Q}\right)$$
$$=\frac{1}{Q^2}\left(1-\frac{2\left(r+1\right)\left(r+2\right)}{r^2}\frac{r}{\left(r+1\right)\left(r+2\right)}\bar{\sigma}^2\psi\right), \tag{47}$$
$$=\frac{1}{Q^2}\left(1-\frac{2}{r}\bar{\sigma}^2\psi\right)$$
$$=\frac{2}{r}\frac{1}{Q^2}\left(\frac{r}{2}-\bar{\sigma}^2\psi\right)$$
$$>0$$

where we used (27). The left inequality in (46) follows from

$$\sqrt{\frac{2\left(r+1\right)}{r+2}\frac{1}{Q}}-\left(\frac{1}{Q}-\frac{\sqrt{1-2Q}}{Q}\right)=\sqrt{\frac{2\left(r+1\right)}{r+2}\frac{1}{Q}}-\frac{1-\sqrt{1-2Q}}{Q}$$
$$=\sqrt{\frac{2\left(r+1\right)}{r+2}\frac{1}{Q}}-\frac{2Q}{Q\left(1+\sqrt{1-2Q}\right)}. \tag{48}$$
$$=\sqrt{\frac{2\left(r+1\right)}{r+2}\frac{1}{Q}}-\frac{2}{1+\sqrt{1-2Q}}$$
$$>\sqrt{\frac{2\left(r+1\right)}{r+2}\frac{1}{Q}}-2$$

Here, we have the following inequality by (27), and the proof is complete:

$$\left(\sqrt{\frac{2(r+1)}{r+2}\frac{1}{Q}}\right)^2-2^2=\frac{2(r+1)}{r+2}\frac{1}{Q}-4$$
$$=\frac{2(r+1)}{r+2}\frac{(r+1)(r+2)}{r\bar{\sigma}^2\psi}-4 \quad (49)$$
$$=4\left\{\frac{(r+1)^2}{2r}\frac{1}{\bar{\sigma}^2\psi}-1\right\}$$
$$>0$$

□

Now, we prove **Proposition 1** stated in the main text.

***Proof of Proposition 1***

The proof is based on the dynamic programming argument (e.g., Section 5 in Baltas and Yannacopoulos [33]). Care is needed regarding the Novikov condition (3) and the global existence of $H$ (more specifically, $A_t$) over time, as shown below. The proof is divided into four steps.

***Step 1. Guess a solution***

Because the problem (10) falls within classical stochastic control problems (e.g., Chapter 5 in Øksendal and Sulem [1]), its optimality equation reads

$$-\frac{\partial H}{\partial t}=\sup_{v\in\mathbb{R}}\left(\left(a_t-\frac{r}{1-t}x+\sigma_t v\sqrt{\frac{r_t}{1-t}x}\right)\frac{\partial H}{\partial x}+\frac{1}{2}\frac{r}{1-t}\sigma_t^2 x\frac{\partial^2 H}{\partial x^2}+x-\frac{1}{2\psi}v^2\right),\quad 0\le t<1,\quad x\ge 0 \quad (50)$$

subject to the terminal condition $H(1,x)=0$, $x\ge 0$. We can rewrite (50) as

$$\begin{aligned}-\frac{\partial H}{\partial t}&=\left(a_t-\frac{r}{1-t}x\right)\frac{\partial H}{\partial x}+\frac{1}{2}\frac{r}{1-t}\sigma_t^2 x\frac{\partial^2 H}{\partial x^2}+x+\sup_{v\in\mathbb{R}}\left(\sigma_t v\sqrt{\frac{r}{1-t}x}\frac{\partial H}{\partial x}-\frac{1}{2\psi}v^2\right)\\&=\left(a_t-\frac{r}{1-t}x\right)\frac{\partial H}{\partial x}+\frac{1}{2}\frac{r}{1-t}\sigma_t^2 x\frac{\partial^2 H}{\partial x^2}+\frac{1}{2}\sigma_t^2\psi\frac{r}{1-t}x\left(\frac{\partial H}{\partial x}\right)^2+x\end{aligned},\quad 0\le t<1,\quad x\ge 0, \quad (51)$$

where

$$v^*(t,x)=\arg\max_{v\in\mathbb{R}}\left(\sigma_t v\sqrt{\frac{r}{1-t}x}\frac{\partial H}{\partial x}-\frac{1}{2\psi}v^2\right)=\sigma_t\psi\sqrt{\frac{r}{1-t}x}\frac{\partial H}{\partial x},\quad 0\le t<1,\quad x\ge 0. \quad (52)$$

Guessing a solution of the form (14) in (51) yields

$$-\frac{\mathrm{d}A_t}{\mathrm{d}t}x-\frac{\mathrm{d}B_t}{\mathrm{d}t}=\left(a_t-\frac{r}{1-t}x\right)A_t+\frac{1}{2}\sigma_t^2\psi\frac{r}{1-t}xA_t^2+x, \quad (53)$$

leading to the following ODEs:

$$-\frac{\mathrm{d}B_t}{\mathrm{d}t}=a_tA_t,\quad 0\le t<1 \quad (54)$$

and

$$-\frac{\mathrm{d}A_t}{\mathrm{d}t}=\frac{r}{1-t}\left(-A_t+\frac{\sigma_t^2\psi}{2}A_t^2\right)+1,\quad 0\le t<1 \tag{55}$$

with terminal conditions $A_1=0$ and $B_1=0$. The ODE (55) is identical to (16) and the relationship (15) follows by integrating (54) considering the terminal condition.

Consequently, we obtain a smooth solution to the optimality equation (50) if three conditions are met: **1) it exists globally in time for the time interval $[0,1]$ (Step 2)**, **2) the sign of $1-\sigma_t^2\psi A_t$ is positive (Step 3)**, and **3) the Novikov condition (3) is satisfied (Step 4)**. By temporally assuming the underlined parts, we can apply the verification result (e.g., Theorem 3.5.2 Pham [26]), showing that the controlled SDE is given by (18) and the optimal control by (17), where the unique existence of the pathwise solution and its nonnegativity is guaranteed by Proposition 1 in Yoshioka [52] due to the specific form of (18). Consequently, what remains to be done is to analyze the **three underlined parts 1)-3)** indicated above, which are addressed in the rest of the proof.

***Step 2. On the nonlinear ODE***

We show that the ODE (55) with the terminal condition $A_1=0$ admits a unique solution that is continuous in $[0,1]$ and continuously differentiable in $(0,1]$. To demonstrate this, we first apply the time transformation $s=1-t$, because dealing with an initial value problem may be more common than dealing with terminal ones. By applying this transformation and setting $U_s=A_{1-t}$, (55) becomes:

$$\frac{\mathrm{d}U_s}{\mathrm{d}s}=\frac{r}{s}\left(-U_s+\frac{\omega_s^2\psi}{2}U_s^2\right)+1,\quad 0<s\le 1 \tag{56}$$

with the initial condition $U_0=0$. Any smooth solution to the initial value problem of the ODE (56) is continuous at the initial time of $s=0$, at which the coefficient on the right-hand side of (56) explodes. Indeed, assuming an asymptotic expansion $U_s\sim U^{(1)}s$ for a small $s>0$ (considering the initial condition $U_0=0$) and substituting it into (56) yields

$$s\frac{\mathrm{d}}{\mathrm{d}s}\left(U^{(1)}s\right)\sim r\left(-U^{(1)}s+\frac{\omega_s^2\psi}{2}\left(U^{(1)}s\right)^2\right)+s\sim -rU^{(1)}s+s\,. \tag{57}$$

Hence, we obtain a first-order expansion of $U_s$ that is valid for a sufficiently small $s>0$:

$$U_s\sim U^{(1)}s=\frac{1}{1+r}s\,. \tag{58}$$

With this result in mind, we apply the following transformation of variables to (56):

$$U_s=\frac{1}{1+r}s\left(1+Y_s\right),\quad 0<s\le 1\,. \tag{59}$$

This leads to

$$\frac{\mathrm{d}Y_s}{\mathrm{d}s}=-\frac{r+1}{s}Y_s+\frac{r}{2(r+1)}\omega_s^2\psi\left(Y_s+1\right)^2,\quad 0<s\le 1 \tag{60}$$

with the initial condition $Y_0 = 0$. Unlike (56), the transformed ODE does not have an exploding coefficient in the nonlinear term. Moreover, if a unique global solution to the initial value problem of transformed ODE (60), then the same applies to (56) due to the asymptotic estimate (58). Therefore, we analyze the ODE (60) in **Step 2**.

Applying a classical variation of the constant formula to (60) along with the initial condition $Y_0 = 0$ yields the following fixed-point problem:

$$\begin{aligned} Y_s &= \frac{r\psi}{2(r+1)} \exp\left( -\int_0^s \frac{r+1}{m} \mathrm{d}m \right) \int_0^s \exp\left( \int_0^m \frac{r+1}{z} \mathrm{d}z \right) \omega_m^2 \left( Y_m + 1 \right)^2 \mathrm{d}m \\ &= \frac{r\psi}{2(r+1)} \frac{1}{s^{r+1}} \int_0^s m^{r+1} \omega_m^2 \left( Y_m + 1 \right)^2 \mathrm{d}m \qquad , \ 0 \le s \le 1 . \\ &\left( \equiv \Theta(Y,s) \right) \end{aligned} \tag{61}$$

We want to show that this fixed-point problem admits a unique solution in the Banach space $C([0,1])$ equipped with the maximum norm $\max_{0\le s\le 1} |\cdot| = |\cdot|_\infty$, but the quadratic nonlinearity $(Y_m + 1)^2$ in $\Theta$ hinders us from directly constructing a contraction mapping because it requires that the target operator, in our case, the integrand of (61), grows at most linearly with respect to $Y_m$. Therefore, we temporarily consider the following relaxed problem and show that it admits a unique nonnegative solution and that the solution is actually the desired solution of (61):

$$Y_s = \frac{r\psi}{2(r+1)} \frac{1}{s^{r+1}} \int_0^s m^{r+1} \omega_m^2 \left( \hat{Y}_m + 1 \right)^2 \mathrm{d}m \ \left( \equiv \hat{\Theta}(Y,s) \right), \ 0 \le s \le 1 \tag{62}$$

with the notation $\hat{Y}_m = \max\left\{0, \min\left\{\bar{Y}, Y_m\right\}\right\}$ for a constant $\bar{Y} > 0$. We choose $\bar{Y}$ so that the four conditions (28), (29), (30), and (31) are satisfied, which is possible by **Lemma 1**.

We demonstrate that the right side of (62) is a strict contraction mapping from $C([0,1])$ to $C([0,1])$. Indeed, for any $Z^{(1)}, Z^{(2)} \in C([0,1])$, we obtain:

$$\begin{aligned} \left| \hat{\Theta}\left( Z^{(1)}, \cdot \right) \right|_\infty &= \max_{0\le s\le 1} \left| \frac{r\psi}{2(r+1)} \frac{1}{s^{r+1}} \int_0^s m^{r+1} \omega_m^2 \left( \hat{Z}_m^{(1)} + 1 \right)^2 \mathrm{d}m \right| \\ &\le \frac{r}{2(r+1)} \bar{\sigma}^2 \psi \left( \bar{Y} + 1 \right)^2 \max_{0\le s\le 1} \left| \frac{1}{s^{r+1}} \int_0^s m^{r+1} \mathrm{d}m \right| \\ &= \frac{r}{2(r+1)} \bar{\sigma}^2 \psi \left( \bar{Y} + 1 \right)^2 \max_{0\le s\le 1} \left| \frac{1}{s^{r+1}} \frac{1}{r+2} s^{r+2} \right| \\ &\le \frac{r}{2(r+1)(r+2)} \bar{\sigma}^2 \psi \left( \bar{Y} + 1 \right)^2 \\ &< +\infty \end{aligned} \tag{63}$$

and

$$
\begin{aligned}
\left|\hat{\Theta}\left(Z^{(1)},\cdot\right)-\hat{\Theta}\left(Z^{(2)},\cdot\right)\right|_{\infty} &= \max_{0\le s\le 1}\left|\frac{r\psi}{2(r+1)}\frac{1}{s^{r+1}}\int_0^s m^{r+1}\omega_m^2\left\{\left(\hat{Z}_m^{(1)}+1\right)^2-\left(\hat{Z}_m^{(2)}+1\right)^2\right\}\mathrm{d}m\right| \\
&\le \frac{r}{2(r+1)}\bar{\sigma}^2\psi\max_{0\le s\le 1}\left(\frac{1}{s^{r+1}}\int_0^s m^{r+1}\left|\left(\hat{Z}_m^{(1)}+1\right)^2-\left(\hat{Z}_m^{(2)}+1\right)^2\right|\mathrm{d}m\right) \\
&\le \frac{r}{2(r+1)}\bar{\sigma}^2\psi\max_{0\le s\le 1}\left(\frac{1}{s^{r+1}}\int_0^s m^{r+1}2\left(\bar{Y}+1\right)\left|\hat{Z}_m^{(1)}-\hat{Z}_m^{(2)}\right|\mathrm{d}m\right) \\
&\le \frac{r}{r+1}\bar{\sigma}^2\psi\left(\bar{Y}+1\right)\max_{0\le s\le 1}\left(\frac{1}{s^{r+1}}\int_0^s m^{r+1}\left|\hat{Z}_m^{(1)}-\hat{Z}_m^{(2)}\right|\mathrm{d}m\right) \\
&\le \frac{r}{r+1}\bar{\sigma}^2\psi\left(\bar{Y}+1\right)\max_{0\le s\le 1}\left(\frac{1}{s^{r+1}}\int_0^s m^{r+1}\left|Z_m^{(1)}-Z_m^{(2)}\right|\mathrm{d}m\right) \\
&\le \frac{r}{r+1}\bar{\sigma}^2\psi\left(\bar{Y}+1\right)\left|Z^{(1)}-Z^{(2)}\right|_{\infty}\max_{0\le s\le 1}\left(\frac{1}{s^{r+1}}\int_0^s m^{r+1}\right) \\
&\le \left(\frac{r}{(r+1)(r+2)}\bar{\sigma}^2\psi\left(\bar{Y}+1\right)\right)\left|Z^{(1)}-Z^{(2)}\right|_{\infty}
\end{aligned}
\quad . \qquad (64)
$$

We have

$$
\begin{aligned}
\left|\hat{\Theta}\left(Z^{(1)},s\right)\right|_{\infty} &\le \frac{r}{2(r+1)}\bar{\sigma}^2\psi\left(\bar{Y}+1\right)^2\left|\frac{1}{s^{r+1}}\frac{1}{r+2}s^{r+2}\right| \\
&\le \frac{r}{2(r+1)(r+2)}\bar{\sigma}^2\psi\left(\bar{Y}+1\right)^2 s \\
&\to 0 \quad (s\to 0)
\end{aligned}
\quad , \qquad (65)
$$

and the continuity of $\hat{\Theta}$ for $0<s\le 1$ follows from its form (62).

According to (64), $\hat{\Theta}$ is a strict contraction mapping because of the first condition (28) in **Lemma 1**. Then, the Banach fixed-point theorem (Theorem 5.7 in Brezis [63]) shows that the relaxed problem (62) admits a unique solution that is bounded, nonnegative, and continuous in $[0,1]$. This solution is referred to as $Z$. Note that the nonnegativity directly follows from the form of $\hat{\Theta}$:

$$
Z_s=\hat{\Theta}(Z,s)=\frac{r\psi}{2(r+1)}\frac{1}{s^{r+1}}\int_0^s m^{r+1}\omega_m^2\left(\hat{Z}_m+1\right)^2\mathrm{d}m\ge 0\ ,\quad 0\le s\le 1\ . \qquad (66)
$$

The upper bound of the solution $Z$ to the relaxed problem is obtained from (62) as done in (63):

$$
|Z|_{\infty}\le\frac{r}{2(r+1)(r+2)}\bar{\sigma}^2\psi\left(\bar{Y}+1\right)^2<\bar{Y}\ , \qquad (67)
$$

where we used the second condition (29) in **Lemma 1**. Consequently, taking "$\max\left\{0,\min\left\{\bar{Y},\cdot\right\}\right\}$" in the relaxed problem is never activated and is hence innocuous. Then, the solution $Z$ to the relaxed problem also solves the original problem (61); hence, we can find a nonnegative and continuous solution to this problem satisfying bound (67). This solution is unique because of the local Lipschitz continuity on the right side of (62). The solution $Z$ is continuously differentiable for $0<s\le 1$ because of

$$\frac{\mathrm{d}}{\mathrm{d}s}\left(\frac{1}{s^{r+1}}\int_0^s m^{r+1}\omega_m^2\left(\hat{Z}_m+1\right)^2\mathrm{d}m\right)=\omega_s^2\left(\hat{Z}_s+1\right)^2-\left(r+1\right)\frac{1}{s^{r+2}}\int_0^s m^{r+1}\omega_m^2\left(\hat{Z}_m+1\right)^2\mathrm{d}m \tag{68}$$

whose right-hand side is continuous in $\left(0,1\right]$. This $Z$ is therefore the desired unique solution.

***Step 3. On the sign of*** $1-\sigma_t^2\psi A_t$

This step is based on **Step 2**. Considering the transformation of variables (59), by the third condition (30) in **Lemma 1**, we have

$$\begin{aligned}\max_{0\le t\le 1}\left|\sigma_t^2\psi A_t\right| &= \max_{0\le s\le 1}\left|\omega_s^2\psi U_s\right| \\ &\le \bar{\sigma}^2\psi\max_{0\le s\le 1}\left|\frac{1}{1+r}s\left(1+Z_s\right)\right| \\ &\le \frac{1}{1+r}\bar{\sigma}^2\psi\max_{0\le s\le 1}\left|\left(1+Z_s\right)\right| \quad , \\ &< \frac{1}{1+r}\bar{\sigma}^2\psi\left(1+\bar{Y}\right) \\ &< 1\end{aligned} \tag{69}$$

and hence

$$\max_{0\le t\le 1}\left|\sigma_t^2\psi A_t\right|<1. \tag{70}$$

From this inequality, the sign of $1-\sigma_t^2\psi A_t$ is entirely positive in $\left[0,1\right]$.

***Step 4. On the Novikov condition***

Finally, we verify the Novikov condition (3). The condition to be checked in the present case is

$$\mathbb{E}_{\mathbb{P}}\left[\exp\left(\frac{1}{2}\int_0^1\left(u_s^*\right)^2\mathrm{d}s\right)\right]<+\infty, \tag{71}$$

which is further rewritten using (17) as follows:

$$\mathbb{E}_{\mathbb{P}}\left[\exp\left(\int_0^1\frac{\sigma_s^2\psi^2A_s^2}{2}\frac{r}{1-s}X_s\mathrm{d}s\right)\right]<+\infty. \tag{72}$$

We set

$$M\left(t,x,k\right)=\mathbb{E}_{\mathbb{P}}\left[\exp\left(K_1\right)\middle|X_t=x,K_t=y\right],\quad 0\le t\le 1,\quad x,k\ge 0 \tag{73}$$

subject to a system of SDEs containing the following auxiliary SDE and the (benchmark) bridge (2):

$$\mathrm{d}K_t=\frac{\sigma_t^2\psi^2A_t^2}{2}\frac{r}{1-t}X_t\mathrm{d}t,\quad 0\le t<1 \tag{74}$$

using the initial conditions of $K_0\ge 0$. Note that the coefficient $\frac{1}{1-t}A_t^2$ is bounded in $\left[0,1\right]$ and vanishes at $t=1$ owing to the asymptotic expansion (58) along with $A_t=U_{1-s}$ and $t=1-s$.

We show that $M\left(0,0,0\right)<+\infty$ and that the Novikov condition is satisfied under the assumption of **Proposition 1**. The Kolmogorov's backward equation associated with (73) is given as follows (note that this equation is seen as an optimality equation without any controls):

$$-\frac{\partial M}{\partial t}=\left(a_t-\frac{r}{1-t}x\right)\frac{\partial M}{\partial x}+\frac{1}{2}\frac{r}{1-t}\sigma_t^2 x\frac{\partial^2 M}{\partial x^2}+\frac{\sigma_t^2\psi^2A_t^2}{2}\frac{r}{1-t}x\frac{\partial M}{\partial y},\quad 0\le t<1,\quad x,k\ge 0 \tag{75}$$

subject to the terminal condition

$$M(1,x,y)=\exp(y),\quad x,y\ge 0. \tag{76}$$

We assume a solution of the following exponential form with smooth time-dependent coefficients $D,E,I$ :

$$M(t,x,y)=\exp\left(D_t x+E_t y+I_t\right),\quad 0\le t\le 1,\quad x,y\ge 0. \tag{77}$$

Substituting (77) into (75) yields

$$-\frac{\mathrm{d}}{\mathrm{d}t}\left(D_t x+E_t y+I_t\right)=\left(a_t-\frac{r}{1-t}x\right)D_t+\frac{1}{2}\frac{r}{1-t}\sigma_t^2 x\left(D_t\right)^2+\frac{\sigma_t^2\psi^2A_t^2}{2}\frac{r}{1-t}xE_t,\quad 0\le t<1,\quad x,y\ge 0. \tag{78}$$

We have the following terminal conditions by comparing (76) and (77):

$$D_1=0,\quad E_1=1,\quad I_1=0. \tag{79}$$

We obtain the following system of ODEs from the equality (78):

$$-\frac{\mathrm{d}D_t}{\mathrm{d}t}=-\frac{r}{1-t}D_t+\frac{1}{2}\frac{r}{1-t}\sigma_t^2\left(D_t\right)^2+\frac{\sigma_t^2\psi^2A_t^2}{2}\frac{r}{1-t}E_t,\quad 0\le t<1, \tag{80}$$

$$-\frac{\mathrm{d}E_t}{\mathrm{d}t}=0,\quad 0\le t<1, \tag{81}$$

$$-\frac{\mathrm{d}I_t}{\mathrm{d}t}=a_tD_t,\quad 0\le t<1. \tag{82}$$

From (79) and (81), we obtain $E_t=1$ ($0\le t\le 1$). From (79) and (82), we obtain:

$$I_t=\int_t^1 a_sD_s\mathrm{d}s,\quad 0\le t\le 1. \tag{83}$$

Because the verification argument is similar to that used in **Step 1**, if $\max_{0\le t\le 1}\left|D_t\right|<+\infty$, then we have $I_0=M(0,0,0)<+\infty$. We show $\max_{0\le t\le 1}\left|D_t\right|<+\infty$ in the rest of this proof.

Because $E_t=1$ ($0\le t\le 1$), the ODE (80) reduces to

$$-\frac{\mathrm{d}D_t}{\mathrm{d}t}=-\frac{r}{1-t}D_t+\frac{1}{2}\frac{r}{1-t}\sigma_t^2\left(D_t\right)^2+\frac{\sigma_t^2\psi^2A_t^2}{2}\frac{r}{1-t},\quad 0\le t<1. \tag{84}$$

Applying the transformation of variables $s=1-t$ and $D_{1-t}=\psi L_s$ to (84) yields

$$\frac{\mathrm{d}\left(\psi L_s\right)}{\mathrm{d}s}=\frac{r}{s}\left(-\psi L_s+\frac{1}{2}\omega_s^2\psi^2L_s^2+\frac{1}{2}\omega_s^2\psi^2U_s^2\right),\quad 0<s\le 1, \tag{85}$$

and hence

$$\begin{aligned}\frac{\mathrm{d}L_s}{\mathrm{d}s}&=\frac{r}{s}\left(-L_s+\frac{1}{2}\omega_s^2\psi L_s^2\right)+\frac{r}{s}\frac{1}{2}\omega_s^2\psi U_s^2\\&=\frac{r}{s}\left(-L_s+\frac{1}{2}\omega_s^2\psi L_s^2\right)+\frac{r}{2(1+r)^2}\omega_s^2\psi\left(1+Y_s\right)^2 s\end{aligned},\quad 0<s\le 1 \tag{86}$$

with the initial condition $L_0=0$, where we used (59). The fourth condition (31) in **Lemma 1** shows

$$\frac{r}{2(1+r)^2}\omega_s^2\psi(1+Y_s)^2 s < \frac{r}{2(1+r)^2}\omega_s^2\psi(1+\bar{Y})^2 s \le s \le 1,\ \ 0 < s \le 1. \tag{87}$$

Any smooth solutions to (85) satisfy $\left.\frac{\mathrm{d}L_s}{\mathrm{d}s}\right|_{s\to+0} = 0$ due to $U_s \sim \frac{s}{r+1}$ for small $s > 0$, and hence $L_s \sim cs^2$ with $c \in \mathbb{R}$ for small $s > 0$. Considering this estimate, applying the transformation of variables

$$L_s = \frac{s}{1+r}W_s,\ \ 0 \le s \le 1 \tag{88}$$

with the initial condition $W_0 = 0$ to (86) yields the following transformed ODE:

$$\frac{\mathrm{d}W_s}{\mathrm{d}s} = -\frac{r+1}{s}W_s + \frac{r}{2(r+1)}\omega_s^2\psi W_s^2 + \frac{r}{2(1+r)}\omega_s^2\psi(1+Y_s)^2,\ \ 0 < s \le 1. \tag{89}$$

From (89), we can consider the fixed-point problem

$$W_s = \frac{r\psi}{2(r+1)}\frac{1}{s^{r+1}}\int_0^s m^{r+1}\omega_m^2 W_m^2 \mathrm{d}m + \frac{r\psi}{2(1+r)}\frac{1}{s^{r+1}}\int_0^s m^{r+1}\omega_m^2(1+Y_m)^2 \mathrm{d}m\ \ \left(=\Xi(W,s)\right),\ \ 0 \le s \le 1, \tag{90}$$

where the second term on the right-hand side of (90) is nonnegative and independent of $W$. Clearly, we have the following *a priori* estimate of any continuous solutions to (90):

$$\min_{0\le s\le 1} W_s \ge 0. \tag{91}$$

We then consider the following auxiliary fixed-point problem:

$$W_s = \tilde{\Xi}(W,s)\left(=\Xi(\tilde{W},s)\right),\ \ 0 \le s \le 1 \tag{92}$$

with the notation $\tilde{W}_m = \max\{0, \min\{\bar{W}, W_m\}\}$ for some $\bar{W} > 0$ chosen later.

Consider any $W^{(1)}, W^{(2)} \in C([0,1])$. Then, as done in **Step 2**, with (87), we have

$$\begin{aligned}\max_{0\le s\le 1}\tilde{\Xi}\left(W^{(1)},\cdot\right) &\le \frac{r\psi}{2(r+1)}\frac{1}{s^{r+1}}\int_0^s m^{r+1}\bar{\sigma}^2\tilde{W}_m^2\mathrm{d}m + (1+r)\frac{1}{s^{r+1}}\int_0^s m^{r+1}\frac{r}{2(1+r)^2}\bar{\sigma}^2\psi(1+\bar{Y})^2\mathrm{d}m \\ &\le \frac{r}{2(r+1)}\bar{\sigma}^2\psi\bar{W}^2\max_{0\le s\le 1}\left|\frac{1}{s^{r+1}}\int_0^s m^{r+1}\mathrm{d}m\right| + (1+r)\max_{0\le s\le 1}\left|\frac{1}{s^{r+1}}\int_0^s m^{r+1}\mathrm{d}m\right| \\ &\le \frac{r}{2(r+1)(r+2)}\bar{\sigma}^2\psi\bar{W}^2 + \frac{r+1}{r+2}\end{aligned} \tag{93}$$

and

$$\min_{0\le s\le 1}\tilde{\Theta}\left(W^{(1)},\cdot\right) \ge 0. \tag{94}$$

Moreover, again, as done in **Step 2**, we have the strict contraction property

$$\left|\tilde{\Xi}\left(W^{(1)},\cdot\right) - \tilde{\Xi}\left(W^{(2)},\cdot\right)\right|_\infty \le \left(\frac{r}{(r+1)(r+2)}\bar{\sigma}^2\psi\bar{W}\right)\left|W^{(1)} - W^{(2)}\right|_\infty \tag{95}$$

for a sufficiently small $\bar{W} > 0$. We can choose $\bar{W} > 0$ such that

$$\frac{r}{2(r+1)(r+2)}\bar{\sigma}^2\psi\bar{W}^2 + \frac{r+1}{r+2} < \bar{W}\ \ \text{and}\ \ \frac{r}{(r+1)(r+2)}\bar{\sigma}^2\psi\bar{W} < 1. \tag{96}$$

In fact, the solution to this inequality system is given by

$$\frac{1}{Q}\left(1-\sqrt{1-2\frac{r+1}{r+2}Q}\right)<\bar{W}<\frac{1}{Q}\text{, where } Q=\frac{r}{(r+1)(r+2)}\bar{\sigma}^2\psi<\frac{1}{2}. \tag{97}$$

The second inequality in (97) comes from **Lemma 1**. In the sequel, we fix one $\bar{W}>0$ such that (97).

Under (96), again by the Banach fixed-point theorem (Theorem 5.7 in Brezis [63]), the auxiliary fixed-point problem (92) admits a unique solution in $C([0,1])$, and moreover, the solution satisfies

$$0\le W_s<\bar{W},\ 0\le s\le 1, \tag{98}$$

showing that the truncation taken in (92) is only surficial. Then, the function $P=W-1$ solves the ODE

$$\frac{\mathrm{d}P_s}{\mathrm{d}s}=-\frac{r+1}{s}P_s+\frac{r}{2(r+1)}\omega_s^2\psi\left(\dot{P}+1\right)^2-\frac{r+1}{s}\left\{1-\frac{r}{2(1+r)^2}\omega_s^2\psi\left(1+Y_s\right)^2 s\right\},\ 0<s\le 1 \tag{99}$$

with the initial condition $P_0=-1$, where we used the notation $\dot{P}_s=\max\left\{-1,\min\left\{\max\left\{\bar{Y},\bar{W}\right\},P_s\right\}\right\}$. The quantity inside $\{\cdot\}$ in (99) is nonnegative due to (87).

Now, the solution $Z$ to (62) also solves the following ODE:

$$\frac{\mathrm{d}Z_s}{\mathrm{d}s}=-\frac{r+1}{s}Z_s+\frac{r}{2(r+1)}\omega_s^2\psi\left(\dot{Z}_s+1\right)^2,\ 0<s\le 1 \tag{100}$$

since

$$-1<0\le Z_s<\bar{Y}\le\max\left\{\bar{Y},\bar{W}\right\},\ 0\le s\le 1. \tag{101}$$

We show that $P_s\le Z_s$ (i.e., $W_s\le Z_s+1$) for $0\le s\le 1$, from which we obtain that $I_0$ is bounded, yielding that the Novikov condition is satisfied. Indeed, if $P_s\le Z_s$ (i.e., $W_s\le Z_s+1$), then

$$\begin{aligned}
I_0&=\int_0^1 a_s D_s\mathrm{d}s\\
&=\int_0^1 a_{1-s}\psi L_s\mathrm{d}s\\
&=\psi\int_0^1 a_{1-s}\frac{s}{1+r}W_s\mathrm{d}s\\
&=\psi\int_0^1 a_{1-s}\frac{s}{1+r}\left(1+P_s\right)\mathrm{d}s\ .\\
&\le\psi\int_0^1 a_{1-s}\frac{s}{1+r}\left(1+Z_s\right)\mathrm{d}s\\
&\le\psi\frac{1+\bar{Y}}{1+r}\left(\max_{0\le s\le 1}a_s\right)\\
&<+\infty
\end{aligned} \tag{102}$$

By (99) and (100), we have

$$\begin{aligned}
\frac{\mathrm{d}\left(P_s-Z_s\right)}{\mathrm{d}s}&\le-\frac{r+1}{s}\left(P_s-Z_s\right)+\frac{r}{2(r+1)}\omega_s^2\psi\left\{\left(\dot{P}_s+1\right)^2-\left(\dot{Z}_s+1\right)^2\right\}\\
&\le-\frac{r+1}{s}\left(P_s-Z_s\right)+\frac{r}{r+1}\bar{\sigma}^2\psi\left(\max\left\{\bar{Y},\bar{W}\right\}+1\right)\left|P_s-Z_s\right|
\end{aligned},\ 0<s\le 1. \tag{103}$$

Assume that there exists some $\bar{s} > 0$ such that $P_s > Z_s$ for $0 < s < \bar{s}$. Then, from (103), we obtain

$$\frac{\mathrm{d}(P_s - Z_s)}{\mathrm{d}s} \leq \left(-\frac{r+1}{s} + \frac{r}{r+1}\bar{\sigma}^2\psi\left(\max\left\{\bar{Y},\bar{W}\right\}+1\right)\right)(P_s - Z_s),\quad 0 < s < \bar{s}, \tag{104}$$

which can be solved for $P_s - Z_s$ to yield the following contradiction:

$$\begin{aligned} P_s - Z_s &\leq (P_0 - Z_0)\exp\left(\int_0^s\left(-\frac{r+1}{m} + \frac{r}{r+1}\bar{\sigma}^2\psi\left(\max\left\{\bar{Y},\bar{W}\right\}+1\right)\right)\mathrm{d}m\right) \\ &= (P_0 - Z_0)s^{-(r+1)}\exp\left(\frac{r}{r+1}\bar{\sigma}^2\psi\left(\max\left\{\bar{Y},\bar{W}\right\}+1\right)s\right) \\ &= -s^{-(r+1)}\exp\left(\frac{r}{r+1}\bar{\sigma}^2\psi\left(\max\left\{\bar{Y},\bar{W}\right\}+1\right)s\right) \\ &< 0 \end{aligned},\quad 0 < s < \bar{s} \tag{105}$$

due to $P_0 = -1$ and $Z_0 = 0$. This implies that such an $\bar{s}$ does not exist, meaning that $P_s \leq Z_s$ for $0 \leq s \leq 1$.

Consequently, we have completed the proof because the **<u>three underlying parts 1)-3)</u>** mentioned in **Step 1** have been resolved.

□